\documentclass[11pt]{amsart}

\usepackage{amsmath}
\usepackage{amsthm}
\usepackage{hyperref}

\newtheorem{prop}{Proposition}
\newtheorem{lemma}[prop]{Lemma}
\newtheorem{thm}[prop]{Theorem}

\theoremstyle{definition}

\newcommand{\dt}{\frac{\partial}{\partial t}}
\newcommand{\brs}[1]{\left| #1 \right|}

\newcommand{\gD}{\Delta}
\newcommand{\gd}{\delta}

\newcommand{\gl}{\lambda}

\newcommand{\ga}{\alpha}

\renewcommand{\ge}{\epsilon}
\newcommand{\N}{\nabla}
\newcommand{\FF}{\mathcal F}

\newcommand{\til}[1]{\widetilde{#1}}
\newcommand{\nm}[2]{|| #1 ||_{#2}}
\newcommand{\gf}{E}

\newcommand{\ohat}[1]{\overset{\circ}{#1}}

\DeclareMathOperator{\Rc}{Rc}
\DeclareMathOperator{\Rm}{Rm}

\DeclareMathOperator{\tr}{tr}

\DeclareMathOperator{\divg}{div}
\DeclareMathOperator{\grad}{grad}

\begin{document}

\title[$L^2$ curvature flow near the round sphere]{The
gradient flow of the $L^2$ curvature energy near the round sphere}
\author{Jeffrey Streets}
\address{Fine Hall\\
         Princeton University\\
         Princeton, NJ 08544}
\email{\href{mailto:jstreets@math.princeton.edu}{jstreets@math.princeton.edu}}

\thanks{The author was partly supported by the National Science Foundation
via
DMS-0703660}

\begin{abstract} We investigate the low-energy behavior of the gradient flow of
the $L^2$ norm of the Riemannian curvature on four-manifolds.  Specifically, we
show long time existence and exponential convergence to a metric of constant
sectional curvature when the initial metric has
positive Yamabe constant and small initial energy.
\end{abstract}

\date{February 21, 2010}

\maketitle

\section{Introduction}

In this paper we study the low-energy behavior of the gradient flow of the
$L^2$ norm of the curvature tensor on four-dimensional manifolds with positive
Yamabe constant.
Let us first introduce some notation.  Let $\Rm$ denote the Riemannian curvature
tensor, $W$ the Weyl curvature, $r$ the Ricci tensor, $z$ the traceless Ricci
tensor, and $s$ the scalar
curvature.  Furthermore, let
\begin{gather*}
\mathcal F(g) := \int_M \brs{\Rm_g}_g^2 dV_g.
\end{gather*}
In what follows we will often drop the explicit reference to $g$, as all objects
in sight will be referencing a given time-dependent metric.  A basic calculation
(\cite{Besse} Proposition 4.70) shows that
\begin{align}
\grad \FF =&\ \gd d r - \check{R} + \frac{1}{4} \brs{\Rm}^2 g.
\end{align}
where $d$ is the exterior derivative acting on the Ricci tensor treated
as a one-form with values in the tangent bundle, and $\gd$ is the adjoint of
$d$.  Moreover, 
\begin{align*}
\check{R}_{ij} = R_{i p q r} R_j^{p q r}.
\end{align*}
A metric is called \emph{critical} if
\begin{align*}
\grad \FF \equiv 0.
\end{align*}
Critical points of quadratic curvature functionals on four-manifolds are very
natural geometric objects to study.  See \cite{LeBrun} for a
nice overview and many interesting results relating the existence of such
metrics to the topology of the underlying manifold.

Given the importance of critical metrics, it is natural to consider the
negative gradient flow of $\mathcal F$:
\begin{gather} \label{flow}
\begin{split}
\dt g =&\ - \grad \FF,\\
g(0) =&\ g_0.
\end{split}
\end{gather}
This is a nonlinear fourth order degenerate parabolic equation.  Since the
equation is fourth order maximum principle techniques are not available,
and the analysis largely relies on integral estimates.  In \cite{Streets} we
showed
short-time existence of the
initial value problem as well as derivative estimates and a long-time existence
obstruction.  Furthermore, in \cite{Streets2} we showed a convergence result for
(\ref{flow}) when the energy is close to zero.  In this paper we examine the
behavior of (\ref{flow}) when the energy is close to its (topologically
determined) minimum and the initial Yamabe constant is positive.  

Before
stating the main result let us fix some further
notation.  Given $M$ a smooth manifold, $\chi(M)$ will denote the Euler
characteristic of $M$.  Also, denote the concircular curvature tensor by
\begin{align*}
\ohat{\Rm} := \Rm - \frac{1}{24} s g \odot g
\end{align*}
where $\odot$ is the Kulkarni-Nomizu product.  Let $(S^4, g_{S^4})$ denote the
sphere with sectional curvature equal to $1$, and likewise $(\mathbb R \mathbb
P^4, g_{\mathbb R \mathbb P^4} )$ is the real projective four-space with
$g_{\mathbb R \mathbb P^4}$ equal to the $\mathbb Z_2$-quotient of $g_{S^4}$. 
Also, for a tensor $T$ we define
\begin{align*}
\brs{\brs{T}}_{L^p} :=&\ \left( \int_M \brs{T}^p \right)^{\frac{1}{p}}, \qquad
\brs{\brs{T}}_{\infty} := \sup_{x \in M} \brs{T}(x)
\end{align*}
and
\begin{align*}
\brs{\brs{T}}^2_{H_k} :=&\ \sum_{j = 0}^k \brs{\brs{\N^j T}}_{L^2}^2.
\end{align*}

\begin{thm} \label{mainthm} There is a constant $\ge > 0$ so that if $(M^4, g)$
is a Riemannian manifold satisfying
\begin{align*}
Y_{[g]} >&\ 0,\\
\brs{\brs{\ohat{\Rm}}}^2_{L^2} \leq&\ \ge \chi(M)
\end{align*}
then then the solution to (\ref{flow}) with initial condition $g$ exists for all
time and converges exponentially to either $(S^4, g_{S^4})$ or $(\mathbb R
\mathbb P^4, g_{\mathbb R \mathbb P^4})$.
\end{thm}

An important remark on the hypotheses is in order.  In dimension $4$, one has
the pointwise equality $\brs{\ohat{\Rm}} = \brs{W}^2 + 2 \brs{z}^2$. 
Therefore the hypothesis includes the
statement that
\begin{align*}
\brs{\brs{W}}_{L^2}^2 \leq \ge \chi(M).
\end{align*}
It follows from \cite{CGY} Theorem A that once $\ge < 16 \pi^2$, $M$ is
diffeomorphic to either $S^4$ or $\mathbb R \mathbb P^4$.  Therefore the theorem
is not providing a new topological conclusion.  Furthermore, the proof relies on
compactness arguments, and so the constant $\ge$ is not computable from the
proof.  There is a discussion of the conjecturally optimal value of $\ge$ in
section 7.  Conversely, the constants $\ge$ in the gap theorems below \emph{are}
computable from the proof, though we do not do this here.

One would hope for an analogous result for metrics with negative scalar
curvature, however the positive sign is crucial for
two main ingredients in the proof.  First of all, in section 2 we exploit a
well-known
relationship between the Yamabe constant, Sobolev constant and the
Gauss-Bonnet formula to show that the hypotheses of a lower bound on the Yamabe
constant and small $L^2$ norm of the traceless curvature tensor imply an
a-priori
estimate of the Sobolev constant.  Next, in section 3, we derive a
coercivity estimate for $\grad \FF$ which holds only for metrics of positive
scalar curvature.  In particular we show that the $H^2$ norm of $\grad \FF$
dominates the $L^2$ norm of $z$.  This estimate is used to show exponential
decay
of the $L^2$ norm of $z$ along solutions to (\ref{flow}), which is one of the
main analytic tools in the proof
of the
theorem.

Another consequence of this key coercivity estimate is the following ``gap
theorem''
for critical metrics, which plays an important role in the proof of Theorem
\ref{mainthm}.  Recall that Einstein metrics and scalar flat,
half-conformally flat metrics are critical for $\mathcal F$.  However, a
complete classification of critical metrics is lacking.  What the following
corollary says is that when the traceless curvature is small in $L^2$ and the
Yamabe constant is positive, a critical metric has constant positive sectional
curvature.

\begin{thm} \label{gapthm} \textbf{Gap Theorem I} There exists $\ge > 0$ 
so that if $(M^4, g)$ is a compact critical Riemannian manifold with
$\brs{\brs{\ohat{\Rm}}}_{L^2}^2 \leq \ge \chi(M)$ and $Y_{[g]} > 0$, then $(M^4,
g)$ is
isometric to $(S^4, g_{S^4})$ or $(\mathbb R \mathbb P^4, g_{\mathbb R \mathbb
P^4})$.
\end{thm}

Related estimates allow us to prove an analogous gap theorem for noncompact
critical manifolds.  This theorem will play a key role in ruling out bubbles in
the proof of Theorem \ref{mainthm}.

\begin{thm} \label{noncompactgapthm} \textbf{Gap Theorem II} Let $(M^4, g)$ be a
noncompact complete
critical Riemannian
four-manifold with zero scalar curvature and $C_S < \infty$.  There is a small
constant $\ge = \ge(C_S) > 0$  so that if $\brs{\brs{\ohat{\Rm}}}^2_{L^2} \leq
\ge$ then $(M^4, g)$ is flat.
\end{thm}

Here is an outline of the rest of the paper.  In section 2 we estimate the
Sobolev constant of metrics with positive Yamabe constant and small $L^2$ norm
of
$\ohat{\Rm}$.  Section 3 contains the main coercive estimate for $\grad \FF$ for
metrics of positive scalar curvature.  Theorem
\ref{gapthm} is a consequence of this estimate, and we finish section 3 by
giving
the proof of Theorem \ref{noncompactgapthm} using related arguments.  In section
4 we give the first
main component of the proof of Theorem \ref{mainthm}, in particular showing that
for $\ge$ chosen small enough solutions to (\ref{flow}) have a definite
lower-bound on their existence time.  This uses an analysis of bubbles,
exploiting Theorem \ref{noncompactgapthm} to rule them out.  After this lower
bound is established one can directly show exponential decay of the energy and
hence convergence of the flow, and this is carried out in sections 5 and 6.  We
conclude in section 7 with some related questions.  Section 8 is an appendix
wherein we show a multiplicative
Sobolev inequality which is used in the proof of the main theorem.

\section{Sobolev Constant Estimate}
In this section we exhibit an estimate of the Sobolev constant of metrics with
positive Yamabe constant and small $L^2$ norm of traceless curvature.  Estimates
of this kind have appeared in many places recently, see \cite{CGY}, \cite{CLW}
for example.  We start by recalling the
Gauss-Bonnet theorem for smooth compact Riemannian four-manifolds:
\begin{align} \label{eq:GB}
\chi(M) =&\ \frac{1}{8 \pi^2} \int_M \left( \frac{s^2}{24} + \brs{W}^2 -
\frac{\brs{z}^2}{2} \right) dV.
\end{align}
Furthermore note that this formula and the conformal invariance
of $\int_M \brs{W}^2$ together imply that
\begin{align*}
\sigma_2(g) := \frac{1}{8 \pi^2} \int_M \left( \frac{s^2}{24} -
\frac{\brs{z}^2}{2}
\right) dV
\end{align*}
is also conformally invariant.

Next recall that the Yamabe constant of a conformal class $[g]$ on a compact
four-manifold is
\begin{align*}
Y_{[g]} = \inf_{\til{g} \in [g]} \frac{\int_M \til{s} d\til{V}}{\left(\int_M d
\til{V} \right)^{\frac{1}{2}}}.
\end{align*}
Applying the solution of the Yamabe problem due to Aubin, Trudinger and Schoen
(\cite{Aubin}, \cite{LeeParker})
this infimum is achieved by a metric of constant scalar curvature.  Using the
expression for the scalar curvature of a conformal metric $\til{g} = u^2 g$ we
conclude
\begin{align*}
Y_{[g]} = \inf_{u \neq 0} \frac{ \int_M \left( 6 \brs{\N u}^2 + s u^2 \right)
dV}{\left( \int_M u^4 dV \right)^{\frac{1}{2}}}.
\end{align*}
In particular it follows that
\begin{align} \label{eq:sob10}
Y_{[g]} \brs{\brs{u}}_{L^4}^2 \leq 6 \brs{\brs{\N u}}_{L^2}^2 + \int_M s u^2 dV
\end{align}
holds for all $u \in C^1(M)$.
Recall that the
\emph{Sobolev
constant} of a metric $g$ on a four-dimensional manifold is the smallest
constant $C_S$ such that the inequality
\begin{align*}
\brs{\brs{u}}_{L^4}^2 \leq C_S \left(\brs{\brs{\N u}}_{L^2}^2 +
V^{-\frac{1}{2}} \brs{\brs{u}}_{L^2}^2 \right)
\end{align*}
holds for all $u \in C^1(M)$.  Now let $(M^4, g)$ be a compact Riemannian manifold
satisfying
\begin{align} \label{eq:sobloc10}
\brs{\brs{\ohat{\Rm}}}_{L^2}^2 \leq \ge \chi(M).
\end{align}
As we noted in the introduction, once $\brs{\brs{\ohat{\Rm}}}_{L^2}^2 \leq 16
\pi^2
\chi(M)$, $M$ is already diffeomorphic to $S^4$ or $\mathbb R \mathbb P^4$.  We
assume for the rest of this section that $M$ is oriented and so $M \cong S^4$,
and so $\chi(M) = 2$.  In particular, using the orthogonal decomposition of the
curvature tensor in dimension $4$, it follows from (\ref{eq:sobloc10}) that
\begin{align*}
\brs{\brs{W}}_{L^2}^2 + 2 \brs{\brs{z}}_{L^2}^2 \leq&\ 2 \ge
\end{align*}
Furthermore it follows from the Gauss-Bonnet theorem that
\begin{align*}
2 + \frac{\ge}{4 \pi^2} \geq\frac{1}{8 \pi^2} \int_M \frac{s^2}{24} \geq&\ 2 -
\frac{\ge}{4 \pi^2}.
\end{align*}
Next it follows from the definition of $\sigma_2$ that
\begin{align*}
\sigma_2(g) > 2 - \frac{\ge}{2 \pi^2}.
\end{align*}
Moreover, due to the conformal invariance of $\sigma_2$, the above inequality holds
for every metric in the conformal class of $g$.  In particular, applying it to
the constant scalar curvature Yamabe minimizer $\til{g}$ we conclude that
\begin{align*}
\frac{1}{192 \pi^2} Y_{[g]}^2 = \frac{1}{192 \pi^2} \frac{\left(\int_M \til{s} d
\til{V} \right)^2}{\int_M d \til{V}} =
\frac{1}{192 \pi^2} \int_M \til{s}^2 d \til{V} \geq&\ \sigma_2(g) > 2 -
\frac{\ge}{2 \pi^2}
\end{align*}
Since $Y_{[g]} > 0$, we conclude $Y_{[g]} > \sqrt{384\pi^2 - 96 \ge}$.  This
also allows us to conclude pinching of the Yamabe energy of $g$.  In particular
we note
\begin{align*}
\int_M s^2 \leq 384 \pi^2 + 96 \ge <&\ Y_{[g]}^2 + 194 \ge\\
\leq&\ \frac{\left(\int_M s dV \right)^2}{\int_M dV} + 194 \ge\\
\leq&\ \int_M s^2 + 194 \ge.
\end{align*}
Let $\bar{s} = \frac{\int_M s dV}{\int_M dV}$. 
Then this estimate implies
\begin{align*}
\brs{\brs{s - \bar{s}}}_{L^2}^2 =&\ \int_M s^2 - \bar{s}^2 \leq 194 \ge.
\end{align*}
Returning to (\ref{eq:sob10}) we may rewrite it as
\begin{align*}
Y_{[g]} \brs{\brs{u}}_{L^4}^2 - \int_M \left(s - \bar{s} \right)u^2 \leq&\ 6
\brs{\brs{\N u}}_{L^2}^2 + \int_M \bar{s} u^2 dV\\
\leq&\ 6 \brs{\brs{\N u}}_{L^2}^2 + V^{-\frac{1}{2}} \sqrt{384 \pi^2 + 48 \ge}
\brs{\brs{u}}_{L^2}^2
\end{align*}
Applying the Cauchy-Schwarz inequality and the above estimate we conclude
\begin{align*}
\int_M \left(s - \bar{s} \right) u^2 \leq&\ \brs{\brs{s - \bar{s}}}_{L^2}
\brs{\brs{u}}_{L^4}^2\\
\leq&\ \sqrt{194 \ge} \brs{\brs{u}}_{L^4}^2.
\end{align*}
Collecting the above estimates together we conclude that if $\ge <
\frac{1}{196}$,
\begin{align*}
\brs{\brs{u}}_{L^4}^2 \leq&\ 768\pi^2 \left( \brs{\brs{\N u}}_{L^2}^2 + V^{-1}
\brs{\brs{u}}_{L^2}^2 \right).
\end{align*}
This completes the proof of the following proposition.
\begin{prop} \label{sobprop} If $(M^4, g)$ is a Riemannian manifold satisfying
\begin{align*}
Y_{[g]} >&\ 0\\
\brs{\brs{\ohat{\Rm}}}_{L^2}^2 \leq&\ \ge \chi(M),
\end{align*}
where $\ge \leq \frac{1}{196}$, then
\begin{align*}
C_S \leq 768 \pi^2
\end{align*}
and
\begin{align*}
Y_{[g]} > \sqrt{384 \pi^2 - 96 \ge}.
\end{align*}
\end{prop}

\section{Coercive Estimate and Gap Theorems}
In this section we exploit the algebraic structure of the gradient of $\mathcal
F$ to derive a coercive estimate for $\grad \FF$ under the assumptions
\begin{gather} \label{coneconds}
\begin{split}
0 < \mu_1 \leq&\ s \leq \mu_2\\
C_S \leq&\ A\\
\brs{\brs{W}}_{L^2}^2 \leq&\ \ge\\
\brs{\brs{z}}_{L^2}^2 \leq&\ \ge\\
\end{split}
\end{gather}
We will apply this estimate to prove Theorem \ref{gapthm}.  Finally we give
the proof of Theorem \ref{noncompactgapthm}, which uses related estimates in the
noncompact setting.  In all estimates below the usage of the constants $\mu_i,
A$, and will
always refer to these constants, while $C$ will denote a generic constant.  The
constant $\ge$ is to be determined by the
estimates
below.  In the end it will depend on $\mu_i$ and $A$ in a way which is
computable in principle, although we do not do this here.

The first step is to derive a partial coercivity estimate from the trace
component of $\grad
\FF$.
\begin{lemma} \label{tracecoercive} There is a constant $C$ so that if
$\epsilon$ is chosen small with
respect to $A$ and $\mu_1$ we have
\begin{align*}
C \brs{\brs{\grad \FF}}_{L^2}^2 \geq\int_M \brs{\N^2 s}^2 + s \brs{\N s}^2.
\end{align*}
\begin{proof} One can directly compute
\begin{align*}
 \tr \grad \FF =&\ - \gD s.
\end{align*}
Therefore
\begin{align*}
4 \brs{\brs{\grad \FF}}_{L^2}^2 \geq&\ \int_M \left(\tr \grad \FF \right)^2\\
=&\ \int_M \N^i \N_i s \N^j \N_j s\\
=&\ - \int_M \N_i s \N^i \N^j \N_j s\\
=&\ - \int_M \N_i s \left( \N^j \N_i \N_j s + R_{j i j}^p \N_p s \right)\\
=&\ \int_M \brs{\N^2 s}^2 + r_{i j} \N_i s \N_j s\\
=&\ \int_M \brs{\N^2 s}^2 + \left( z_{i j} + \frac{1}{4} s g_{i j} \right) \N^i
s \N^j s\\
=&\ \int_M \brs{\N^2 s}^2 + z_{i j} \N^i s \N^j s + \frac{1}{4} s \brs{\N s}^2
\end{align*}
Next by applying H\"older's inequality we estimate
\begin{align*}
 \int_M z_{i j} \N^i s \N^j s \leq&\ \brs{\brs{z}}_{L^2} \brs{\brs{ \brs{\N
s}^2}}_{L^2}\\
\leq&\ \ge^{\frac{1}{2}} \left(\int_M \brs{\N s}^4 \right)^{\frac{1}{2}}\\
\leq&\ A \ge^{\frac{1}{2}} \left( \int_M \brs{\N^2 s}^2 + \int_M \brs{\N
s}^2 \right)\\
\leq&\ A \ge^{\frac{1}{2}} \left( \int_M \brs{\N^2 s}^2 + \frac{1}{\mu_1} \int_M
s \brs{\N s}^2 \right).
\end{align*}
Thus for $\ge$ chosen small enough with respect to $A$ and $\mu_1$, the result
follows.
\end{proof}
\end{lemma}

Next we derive a coercivity estimate from the full tensor $\grad \FF$.  Before
the proof we will record a special expression for $\grad \FF$ in four
dimensions.
\begin{lemma} \label{fullcoercivelemma2} Let $(M^4, g)$ be a Riemannian
four-manifold.  Then
\begin{align*}
 \grad \FF =&\ -2 \gD r + \N^2 s + \frac{s}{3} z + 4 z \circ z - \brs{z}^2 g - 4
W \circ z.
\end{align*}
\begin{proof}  First of all (\cite{Besse} Proposition 4.70) implies that
\begin{align*}
\grad \FF = \gd d r - \check{R} + \frac{1}{4} \brs{\Rm}^2 g.
\end{align*}
Next note that, in four dimensions (\cite{Besse} (4.72)),
\begin{align*}
\check{R} - \frac{1}{4} \brs{\Rm}^2 g = \frac{s}{3} z + 2 W \circ z.
\end{align*}
Also,
\begin{align*}
\gd d r =&\ - 2 \gD r + \N^2 s + 2 r \circ r - 2 R \circ r.
\end{align*}
Combining these yields
\begin{align*}
\grad \FF =&\ - 2 \gD r + \N^2 s + 2 r \circ r - 2 R \circ r - \frac{s}{3} z - 2
W \circ z
\end{align*}
Now we write
\begin{align*}
2 r \circ r =&\ 2 \left(z + \frac{s}{4} g \right) \circ \left(z + \frac{s}{4} g
\right)\\
=&\ 2 z \circ z + s z + \frac{1}{8} s^2 g.
\end{align*}
Also, recall the four-dimensional curvature decomposition
\begin{gather} \label{eq:curvdec}
\begin{split}
R_{i j k l} =&\ W_{i j k l} + \frac{1}{2} \left( z_{i l} g_{j k} - z_{i k} g_{j
l} + z_{j k} g_{i l} - z_{j l} g_{i k} \right)\\
&\ + \frac{1}{12} s \left( g_{j k} g_{i l} - g_{j l} g_{i k} \right).
\end{split}
\end{gather}
We conclude that
\begin{align*}
 - 2 R \circ r =&\ - 2 W \circ z - \brs{z}^2 g + 2 z \circ z - \frac{s}{3} z -
\frac{1}{8} s^2 g.
\end{align*}
Combining these calculations yields the result.
\end{proof} 
\end{lemma}

\begin{prop} \label{maincoerciveestimate} Given $\mu_i, A > 0$, there are
constants $\gd = \gd(\mu_i)
> 0$ and $\ge = \ge(\mu_i, A)$ such that if $(M^4, g)$ is a
compact Riemannian manifold satisfying (\ref{coneconds}), then we have
\begin{align*}
\brs{\brs{\grad \FF}}_{L^2}^2 \geq \gd \left( \brs{\brs{\gD r}}_{L^2}^2 +
\brs{\brs{z}}_{H^1}^2 \right).
\end{align*}
\begin{proof} We start with the result of Lemma \ref{fullcoercivelemma2} and
expand the $L^2$ inner product
\begin{align*}
 \brs{\brs{\grad \FF}}_{L^2}^2 =&\ \brs{\brs{ -2 \gD r + \N^2 s + \frac{s}{3} z
+ 4 z
\circ z - \brs{z}^2 g - 4 W \circ z}}_{L^2}^2\\
=&\ 4 \brs{\brs{\gD r}}_{L^2}^2 + \brs{\brs{\N^2 s}}_{L^2}^2 +
\brs{\brs{\frac{s}{3} z}}_{L^2}^2\\
&\ + \brs{\brs{4 z \circ z - \brs{z}^2 g}}_{L^2}^2 + 16 \brs{\brs{W \circ
z}}_{L^2}^2\\
&\ - 4 \left< \gD r, \N^2 s \right>_{L^2} - 4 \left<\gD r, \frac{s}{3} z
\right>_{L^2} - 4
\left< \gD r, 4 z \circ z - \brs{z}^2 g \right>_{L^2}\\
&\ + 16 \left< \gD r, W \circ z \right> + 2 \left< \N^2 s, \frac{s}{3} z
\right>_{L^2} + 2 \left< \N^2 s, 4 z \circ z -
\brs{z}^2 g \right>_{L^2}\\
&\ - 8 \left<\N^2 s, W \circ z \right> + 2 \left< \frac{s}{3} z, 4 z \circ z
\right>_{L^2} - 8 \left< \frac{s}{3} z, W \circ z \right>\\
=:&\ \sum_{j = 1}^{14} I_j.
\end{align*}
We now estimate the individual terms $I_j$.
First using Lemma \ref{tracecoercive} we conclude that
\begin{align*}
I_6 \geq&\ - \theta \brs{\brs{\gD r}}_{L^2}^2 - \frac{C}{\theta} \brs{\brs{\grad
\FF}}_{L^2}^2
\end{align*}
where $\theta$ is a small constant to be determined later.
Next consider
\begin{align*}
 I_7 =&\ - \frac{4}{3} \int_M \left< \gD r, s z \right>\\
=&\ \frac{4}{3} \left< \N r, z \N s + s \N z \right>\\
=&\ \frac{4}{3} \int_M s \brs{\N z}^2 + \int_M \N z * \N s * z.
\end{align*}
Now using that $\N s$ may be expressed in terms of $\N z$ by the
Bianchi identity we estimate
\begin{gather} \label{eq:coerloc2}
\begin{split}
\brs{\int_M \N z * \N s * z} \leq&\ C \brs{\brs{z}}_{L^2} \brs{\brs{ \brs{\N
z}^2}}_{L^2}\\
\leq&\ C \ge^{\frac{1}{2}} \left(\int_M \brs{\N z}^4 \right)^{\frac{1}{2}}\\
\leq&\ C A \ge^{\frac{1}{2}} \left( \int_M \brs{\N^2 z}^2 + \int_M \brs{\N
z}^2 \right).
\end{split}
\end{gather}
Thus for $\ge$ chosen small with respect to $A$ we conclude
\begin{align*}
 I_7 \geq&\ \frac{4}{3} \int_M s \brs{\N
z}^2 - C A \ge^{\frac{1}{2}} \brs{\brs{z}}_{H^2}^2.
\end{align*}
Similar estimates yield
\begin{align*}
I_8 + I_{11} \geq - C A \ge^{\frac{1}{2}} \brs{\brs{z}}_{H^2}^2.
\end{align*}
Next consider
\begin{align*}
I_9 \geq&\ - \theta \brs{\brs{\gD r}}_{L^2}^2 - \frac{C}{\theta} \brs{\brs{W
\circ z}}_{L^2}^2.
\end{align*}
Next we estimate, again using the Bianchi identity,
\begin{align*}
\brs{I_{10}} =&\ 2 \brs{\int_M \left< \N^2 s, \frac{s}{3} z \right>}\\
=&\ \frac{2}{3} \brs{\int_M \N^i \N^j s  \left( s z_{i j} \right)}\\
=&\ \brs{\int_M z * \N s^{*2} + s * \N s^{*2}}\\
\leq&\ C \int_M \brs{z} \brs{\N s}^2 + s \brs{\N s}^2.
\end{align*}
From (\ref{eq:coerloc2}) and Lemma \ref{tracecoercive} we conclude
\begin{align*}
I_{10} \geq&\ - C A \ge^{\frac{1}{2}} \brs{\brs{z}}_{H^2}^2 - C \brs{\brs{\grad
\FF}}_{L^2}^2.
\end{align*}
Next applying Lemma \ref{tracecoercive} we conclude
\begin{align*}
I_{12} \geq&\ - \theta \brs{\brs{W \circ z}}_{L^2}^2 - \frac{C}{\theta}
\brs{\brs{\grad \FF}}_{L^2}^2
\end{align*}
Next we have
\begin{align*}
\brs{I_{13}} =&\ \brs{\int_M s z^{*3}}\\
\leq&\ C \brs{\brs{s z}}_{L^2} \left( \int_M \brs{z}^4 \right)^{\frac{1}{2}}\\
\leq&\ C A \brs{\brs{s z}}_{L^2} \left( \int_M \brs{\N z}^2 + \int_M \brs{z}^2
\right).
\end{align*}
Next we estimate
\begin{align*}
C A \brs{\brs{s z}}_{L^2} \brs{\brs{z}}_{L^2}^2 =&\ C A \ge^{\frac{1}{2}} \left(
\brs{\brs{ s z}}_{L^2} \brs{\brs{z}}_{L^2} \right)\\
\leq&\ C A \ge^{\frac{1}{2}} \left( \brs{\brs{ s z}}_{L^2}^2 +
\brs{\brs{z}}_{L^2}^2 \right).
\end{align*}
Also,
\begin{align*}
C A\brs{\brs{s z}}_{L^2} \int_M \brs{\N z}^2 =&\ - C A \brs{\brs{s z}}_{L^2}
\int_M
\left<
z, \N^2 z \right>\\
\leq&\ C A \brs{\brs{s z}}_{L^2} \brs{\brs{z}}_{L^2} \brs{\brs{\N^2 z}}_{L^2}\\
\leq&\ C A \ge^{\frac{1}{2}} \left( \brs{\brs{s z}}_{L^2}^2 +
\brs{\brs{z}}_{H^2}^2 \right).
\end{align*}
Combining these we conclude
\begin{align*}
I_{13} \geq&\ - C A \ge^{\frac{1}{2}} \left(\brs{\brs{s z}}_{L^2}^2 +
\brs{\brs{z}}_{H^2}^2
\right).
\end{align*}
Finally we estimate
\begin{align*}
I_{14} \geq&\ - \theta \brs{\brs{s z}}_{L^2}^2 - \frac{C}{\theta} \brs{\brs{W
\circ z}}_{L^2}^2.
\end{align*}
Collecting these individual estimates, choosing $\ge$ small with respect to $A$,
and choosing $\theta$ small yields
\begin{gather} \label{eq:coerloc3}
\begin{split}
C \brs{\brs{\grad \FF}}_{L^2}^2 \geq&\ \frac{1}{10} \left( \brs{\brs{\gD
r}}_{L^2}^2 +
\brs{\brs{\N^2 s}}_{L^2}^2 + \brs{\brs{s z}}_{L^2}^2 + \int_M s \brs{\N z}^2
\right)\\
&\ - C A \ge^{\frac{1}{2}} \brs{\brs{z}}_{H^2}^2 - C
\brs{\brs{W \circ z}}_{L^2}^2.
\end{split}
\end{gather}
It remains to estimate the last two terms, which we do in the next two lemmas. 
First consider
\begin{lemma} \label{le:Weyllemma2} We may choose $\ge$ small with respect to
$A$ so that
\begin{align*}
 \brs{\brs{\brs{W}\brs{z}}}_{L^2}^2 \leq&\ C(\mu_i, A) \ge \brs{\brs{z}}_{H^2}^2
\end{align*}
\begin{proof}
Begin by applying H\"older's inequality and the Sobolev inequality to yield
\begin{gather} \label{eq:Weyllemmaloc5}
\begin{split}
\int_M \brs{W}^2 \brs{z}^2 \leq&\ \left( \int_M \brs{W}^4 \right)^{\frac{1}{2}}
\left( \int_M \brs{z}^4 \right)^{\frac{1}{2}}\\
\leq&\ A^2 \left( \int_M \brs{\N W}^2 + \int_M \brs{W}^2 \right) \left( \int_M
\brs{\N z}^2 + \int_M \brs{z}^2 \right)\\
=&\ A^2 \left( \brs{\brs{\N W}}_{L^2}^2 + \brs{\brs{W}}_{L^2}^2 \right) \left(
\brs{\brs{z}}_{H^1}^2 \right).
\end{split}
\end{gather}
Before continuing we need a nice expression for $\gD W$.  First apply the second
Bianchi identity and commute derivatives to yield
\begin{gather} \label{eq:coerloc15}
\begin{split}
\N_i \N_i W_{m k l n} =&\ \N_i \left(\N_m W_{i k l n} + \N_k W_{m i l n}
\right) + \N^2 z\\
=&\ \N_m \N_i W_{i k l n} + \N_k \N_i W_{m i l n} + \Rm * W + \N^2 z\\
=&\ \N_m \N_k W_{i i l n} + \N_k \N_m W_{i i l n} + \Rm * W + \N^2 z\\
=&\ \N^2 z + \Rm * W
\end{split}
\end{gather}
We conclude
\begin{align*}
\brs{\brs{\N W}}_{L^2}^2 =&\ - \int_M \left<W, \gD W \right>\\
=&\ \int_M \Rm * W^{*2} + \int_M W * \N^2 z.
\end{align*}
To estimate the first term in the line above, note
\begin{align*}
\int_M s * W^{*2} \leq&\ C \mu_2 \brs{\brs{W}}_{L^2}^2,
\end{align*}
Also we have
\begin{align*}
\int_M W^{*3} \leq&\ C \brs{\brs{W}}_{L^2} \left(\int_M \brs{W}^4
\right)^{\frac{1}{2}}\\
\leq&\ C A \ge^{\frac{1}{2}} \left( \brs{\brs{\N W}}_{L^2}^2 +
\brs{\brs{W}}_{L^2}^2 \right).
\end{align*}
There is also the curvature term
\begin{align*}
\int_M z * W^{*2} \leq&\ C \brs{\brs{z}}_{L^2} \left( \int_M \brs{W}^4
\right)^{\frac{1}{2}}\\
\leq&\ C A \ge^{\frac{1}{2}} \left( \brs{\brs{\N W}}_{L^2}^2 +
\brs{\brs{W}}_{L^2}^2 \right).
\end{align*}
Collecting these calculations and choosing $\ge$ small with respect to $A$ we
conclude
\begin{align*}
\brs{\brs{\N W}}_{L^2}^2 \leq&\ C(\mu_2) \brs{\brs{W}}_{L^2}^2
+ \int_M W * \N^2 z.
\end{align*}
Plugging this into (\ref{eq:Weyllemmaloc5}) yields
\begin{align*}
 \int_M \brs{W}^2 \brs{z}^2 \leq&\ C(\mu_2) A^2 \ge \brs{\brs{z}}_{H^1}^2 + A^2
\left(\int_M W * \N^2 z \right) \brs{\brs{z}}_{H_1}^2.
\end{align*}
To estimate the final term we first consider
\begin{align*}
 A^2 \left( \int_M W * \N^2 z \right) \brs{\brs{\N z}}_{L^2}^2 \leq&\ C A^2
\brs{\brs{W}}_{L^2} \brs{\brs{\N^2 z}}_{L^2} \brs{\brs{\N z}}_{L^2}^2\\
\leq&\ C A^2 \brs{\brs{W}}_{L^2} \brs{\brs{z}}_{L^2} \brs{\brs{\N^2
z}}_{L^2}^2\\
\leq&\ C A^2 \ge \brs{\brs{\N^2 z}}_{L^2}^2.
\end{align*}
Lastly estimate
\begin{align*}
A^2 \left( \int_M W * \N^2 z \right) \brs{\brs{z}}_{L^2}^2 \leq&\ C A^2
\brs{\brs{W}}_{L^2} \brs{\brs{\N^2 z}}_{L^2} \brs{\brs{z}}_{L^2}^2\\
\leq&\ C A^2 \ge \brs{\brs{\N^2 z}}_{L^2}
\brs{\brs{z}}_{L^2}\\
\leq&\ C A^2 \ge \left( \brs{\brs{\N^2
z}}_{L^2}^2 + \brs{\brs{z}}_{L^2}^2 \right).
\end{align*}
The result follows.
\end{proof}
\end{lemma}

\begin{lemma} \label{fullcoercivelemma1} There is a constant $C$ so that if we
choose $\ge$ small with respect
to $A$,
\begin{align*}
\brs{\brs{\N^2 z}}_{L^2}^2 \leq&\ C \left( \brs{\brs{\gD z}}_{L^2}^2 +
\brs{\brs{z}}_{H^1}^2 + \brs{\brs{s z}}_{L^2}^2 + \int_M s \brs{\N z}^2 +
\brs{\brs{\brs{W}\brs{z}}}^2_{L^2} \right).
 \end{align*}
\begin{proof}
 We integrate by parts and estimate
\begin{align*}
\brs{\brs{\N^2 z}}_{L^2}^2 =&\ \int_M \N_i \N_j z_{kl} \N_i \N_j z_{kl}\\
=&\ \int_M \N_j \N_i z_{kl} \N_i \N_j z_{k l} + \Rm * z * \N^2 z\\
\leq&\ - \int_M \N_i z_{k l} \N_j \N_i \N_j z_{k l} + C \int_M \brs{\Rm}^2
\brs{z}^2 +
\frac{1}{2} \int_M \brs{\N^2 z}^2\\
=&\ - \int_M \N_i z_{k l} \N_i \gD z_{k l} + \int_M \Rm * \N z^{*2} + C \int_M
\brs{\Rm}^2 \brs{z}^2 + \frac{1}{2} \int_M \brs{\N^2 z}^2\\
\leq&\ \brs{\brs{\gD z}}_{L^2}^2 + C \int_M \brs{\Rm} \brs{\N z}^2 +
C \int_M \brs{\Rm}^2 \brs{z}^2 + \frac{1}{2} \int_M
\brs{\N^2 z}^2\\
\leq&\ C \left( \brs{\brs{\gD z}}_{L^2}^2 + \int_M \brs{\Rm} \brs{\N z}^2 +
 \int_M \brs{\Rm}^2 \brs{z}^2 \right).
\end{align*}
Next we estimate
\begin{align*}
\int_M \brs{\Rm} \brs{\N z}^2 \leq&\ \int_M s \brs{\N z}^2 + \int_M \left(
\brs{W} + \brs{z} \right) \brs{\N z}^2\\
\leq&\ \int_M s \brs{\N z}^2 + \left( \brs{\brs{W}}_{L^2} +
\brs{\brs{z}}_{L^2} \right) \left( \int_M \brs{\N z}^4 \right)^{\frac{1}{2}}\\
\leq&\ \int_M s \brs{\N z}^2 + A \left( \brs{\brs{W}}_{L^2} +
\brs{\brs{z}}_{L^2} \right) \left( \int_M \brs{\N^2 z}^2 + \int_M \brs{\N z}^2
\right)\\
\leq&\ \int_M s \brs{\N z}^2 + A \ge^{\frac{1}{2}} \left( \brs{\brs{\N^2
z}}_{L^2}^2 + \brs{\brs{\N z}}_{L^2}^2 \right)
\end{align*}
Also we estimate
\begin{align*}
\int_M \brs{\Rm}^2 \brs{z}^2 \leq&\ \brs{\brs{s z}}_{L^2}^2 + \int_M \brs{z}^4 +
\brs{\brs{ \brs{W} \brs{z}}}_{L^2}^2\\
\leq&\ \brs{\brs{s z}}_{L^2}^2 + 2 \left(\int_M \brs{\N z}^2 \right)^2 + 2
\left( \int_M \brs{z}^2 \right)^2 + \brs{\brs{ \brs{W} \brs{z}}}_{L^2}^2\\
\leq&\ \brs{\brs{s z}}_{L^2}^2 + 2 \ge \brs{\brs{\N^2 z}}_{L^2}^2 + 2 \ge
\brs{\brs{z}}_{L^2}^2 + \brs{\brs{ \brs{W} \brs{z}}}_{L^2}^2.
\end{align*}
Combining these estimates yields the result.
\end{proof}
\end{lemma}
Applying Lemmas \ref{le:Weyllemma2} and \ref{fullcoercivelemma1}, and the fact
that $s > \mu_1$, we conclude from (\ref{eq:coerloc3}) that if $\ge$ is chosen
small we have
\begin{align*}
C \brs{\brs{ \grad \FF}}_{L^2}^2 \geq&\ \frac{1}{20} \brs{\brs{\gD r}}_{L^2}^2 +
\left( \frac{\mu_1}{20} - C(A, \mu_i) \ge \right) \brs{\brs{\N z}}_{L^2}^2\\
&\ + \left( \frac{\mu^2_1}{20} - C(A, \mu_i) \ge \right)
\brs{\brs{z}}_{L^2}^2.
\end{align*}
The proposition follows.
\end{proof}
\end{prop}

\begin{thm} \textbf{Gap Theorem I} There exists $\ge > 0$ 
so that if $(M^4, g)$ is a compact critical Riemannian manifold with
$\brs{\brs{\ohat{\Rm}}}_{L^2}^2 \leq \ge \chi(M)$ and $Y_{[g]} > 0$, then $(M^4,
g)$ is
isometric to $(S^4, g_{S^4})$ or $(\mathbb R \mathbb P^4, g_{\mathbb R \mathbb
P^4})$.
\begin{proof}  As pointed out in the introduction, the hypotheses already imply
that $M$ is diffeomorphic to $S^4$ or $\mathbb R \mathbb P^4$, and by passing to
the double cover we may as well assume $M \cong S^4$.  Since $\tr_g \grad \FF =
\gD s$ and $g$ is critical it follows that $g$ has constant
scalar curvature.  Scale $g$ so that it has unit volume, then we have $s \equiv
Y_{[g]}$.  Apply Proposition \ref{sobprop} to conclude that the Sobolev constant
of $g$ is bounded above and the Yamabe constant is bounded below, so the scalar
curvature is bounded above and below.  We now may apply Proposition
\ref{maincoerciveestimate} to conclude that if $\brs{\brs{\ohat{\Rm}}}_{L^2}$ is
chosen small enough with respect to $Y_{[g]}$ then $g$ is Einstein.  It follows
from the arguments of section 2 that in fact for $\ge$ small the Yamabe constant
of $(M^4, g)$ is close to that of $S^4$.  It now follows from Theorem C of
\cite{Gursky} that $(M^4, g)$ is isometric to $(S^4, g_{S^4})$.  We also sketch 
another argument below to finish the theorem which is more in line with the type
of arguments we have been using.

Since $g$ is now Einstein, it follows
that the traceless part of the curvature tensor satisfies the elliptic equation
\begin{align*}
\gD \ohat{\Rm} =&\ \ohat{\Rm} * \ohat{\Rm} + s * \ohat{\Rm}.
\end{align*}
It follows that the curvature satisfies the local elliptic estimate
\begin{align*}
\sup_{B_r} \brs{\ohat{\Rm}} \leq \frac{C}{r} \brs{\brs{\ohat{\Rm}}}_{L^2(B_r)}
\end{align*}
for balls satisfying
\begin{align*}
 \brs{\brs{\ohat{\Rm}}}_{L^2(B_r)} \leq \ge_0.
\end{align*}
The constants $\ge_0$ and $C$ here depend on a bound for $s$ and a bound on the
Sobolev constant, both of which are bounded by $Y_{[g]}$.  Therefore, for $\ge$
chosen small enough, $\brs{\brs{\ohat{\Rm}}}_{\infty} \leq
C(Y_{[g]}) \ge$.  In particular, for $\ge$ chosen small enough we can conclude
that $g$ has positive curvature operator.  It now follows from the main theorem
of \cite{Hamilton4} that in fact $g$ is isometric to the round metric on $S^4$.
\end{proof}
\end{thm}

\begin{thm} \textbf{Gap Theorem II} Let $(M^4, g)$ be a noncompact complete
critical Riemannian
four-manifold with zero scalar curvature and $C_S < \infty$.  There is a small
constant $\ge = \ge(C_S) > 0$  so that if $\brs{\brs{\ohat{\Rm}}}^2_{L^2} \leq
\ge$ then $(M^4, g)$ is flat.
\begin{proof}
Since $s \equiv 0$, let us write the critical equation in the simple form
\begin{align*}
0 =&\ \gD r + \Rm^{*2}
\end{align*}
Let $\phi$ be some compactly supported function.  First observe the inequality
\begin{align*}
- \brs{\Rm} \gD \brs{\Rm} =&\ - \frac{1}{2} \gD \brs{\Rm}^2 + \brs{\N \brs{\Rm}}\\
=&\ - \left< \gD \Rm, \Rm \right> - \brs{\N \Rm}^2 + \brs{\N \brs{\Rm}}\\
\leq&\ - \left< \gD \Rm, \Rm \right>
\end{align*}
using the Kato inequality $\brs{\N \brs{\Rm}} \leq \brs{\N \Rm}$.  By the Bianchi identity, one can show that $\gD \Rm = \mathcal L(\N \divg \Rm) + \Rm * \Rm$ for some universal linear operator $\mathcal L$.  Therefore we may estimate for any Riemannian metric
\begin{align*}
- \int_M \phi^2 \brs{\Rm} \gD \brs{\Rm} dV \leq&\ \int_M \phi^2 \left( - \left<\gD \Rm, \Rm \right> \right) dV\\
=&\ \int_M \phi^2 \left(- \left<\mathcal  L(\N \divg \Rm) + \Rm^{*2}, \Rm \right> \right) dV\\
=&\ \int_M \left(\phi * \N \phi * \divg \Rm * \Rm + \phi^2 \divg \Rm^{*2} + \phi^2 * \Rm^{*3} \right) dV\\
\leq&\ C \int_M \brs{\N \phi}^2
\brs{\Rm}^2 + \phi^2 \brs{\divg \Rm}^2 + \brs{\Rm}^3 \phi^2 dV.
\end{align*}
Next we use the critical equation to estimate
\begin{align*}
0 =&\ \int_M \phi^2 \left< \gD r + \Rm^{*2}, r \right> dV\\
=&\ - \int_M \phi^2 \brs{\N r}^2 + \phi r * \N r * \N \phi + \phi^2 \Rm^{*3}
dV\\
\leq&\ - \frac{1}{2} \int_M \phi^2 \brs{\N r}^2 dV + C \int_M \left(\brs{\N
\phi}^2 \brs{\Rm}^2 + \phi^2 \brs{\Rm}^3 \right)dV.
\end{align*}
Since $\divg \Rm = \N^i \Rm_{i j k l} = \N_k r_{j l} - \N_l r_{j k}$ by the
Bianchi identity, we conclude that
\begin{align*}
- \int_M \phi^2 \brs{\Rm} \gD \brs{\Rm} dV \leq&\ C \int_M \left(\brs{\N \phi}^2
\brs{\Rm}^2 + \phi^2 \brs{\Rm}^3 \right) dV.
\end{align*}
Applying the Sobolev inequality we conclude using the above estimate that
\begin{align*}
\brs{\brs{ \phi \brs{\Rm}}}_{L^4}^2 \leq&\ C \int_M \left(\brs{\N \phi}^2
\brs{\Rm}^2 + \brs{\N \brs{\Rm}}^2 \phi^2 \right) dV\\
\leq&\ C \int_M \left(\brs{\N \phi}^2 \brs{\Rm}^2 + \phi^2 \brs{\Rm}^3
\right) dV\\
\leq&\ C \int_M \brs{\N \phi}^2 \brs{\Rm}^2 dV + C \brs{\brs{ \phi
\brs{\Rm}}}_{L^4}^2 \brs{\brs{\Rm}}_{L^2}.
\end{align*}
Therefore for $\ge$ chosen small with respect to the Sobolev constant we
conclude
\begin{align*}
\brs{\brs{\phi \brs{\Rm}}}_{L^4}^2 \leq C \int_M \brs{\N \phi}^2 \brs{\Rm}^2
dV.
\end{align*}
Fix some point $x \in M$, and let $\phi$ be a cutoff function for the ball of
radius $\rho$.  In particular choose $\phi$ such that
\begin{align*}
0 \leq \phi \leq&\ 1\\
\phi =&\ 1 \mbox{ on } B_{\frac{\rho}{2}}(x)\\
\phi =&\ 0 \mbox{ on } M \setminus B_{\rho}(x)\\
\brs{\N \phi} \leq&\ \frac{4}{\rho}.
\end{align*}
It follows that
\begin{align*}
\brs{\brs{\phi \brs{\Rm}}}_{L^4}^2 \leq&\ \frac{4}{\rho^2} \int_{B_{\rho} -
B_{\frac{\rho}{2}}} \brs{\Rm}^2\\
\leq&\ \frac{C}{\rho^2}.
\end{align*}
Letting $\rho \to \infty$ we conclude $\brs{\Rm} \equiv 0$, and the result
follows.
\end{proof}
\end{thm}

\section{Proof of Theorem 1}

\begin{proof}
We proceed by contradiction.  If the statement is false, then we may choose
$\ge_i \to 0$ and metrics $g_i$ such that $\brs{\brs{\ohat{\Rm_{g_i}}}}_{L^2}^2
\leq
\ge_i \chi(M)$, and the solution to (\ref{flow}) with initial condition $g_i$
exists on
a finite time interval.  As noted in the introduction, once $\ge_i < 16 \pi^2$
it follows that $M$ is diffeomorphic to either $S^4$ or $\mathbb R \mathbb
P^4$, so we can conclude that $\chi(M) = 2, 1$.  By lifting to the double cover,
we may assume without loss of generality that $M \cong S^4$, and by redefining
$\ge_i$ that
\begin{align*}
\brs{\brs{\ohat{\Rm_{g_i}}}}_{L^2} \leq \ge_i.
\end{align*}

The first major step is to use a blowup argument to
show that the existence time is bounded below for $i$ sufficiently large.  It is
important to note here that the small energy condition above is not a priori
preserved in general for solutions to (\ref{flow}).  Indeed, it follows from the
Gauss-Bonnet theorem (\cite{Ber}) that
\begin{align} \label{eq:energy}
\mathcal F(g) =&\ 8 \pi^2 \chi(M) + \int_M \brs{z}^2 dV.
\end{align}
Therefore an upper bound on $\brs{\brs{z}}_{L^2}^2$ is automatically preserved,
but it is possible that the balance between the scalar curvature and Weyl
curvature contributions to $\mathcal F$ could change along the flow.  This
important technical issue
is discussed in some more detail in section 7.  To control the balance between
scalar and Weyl curvatures we
first need a lemma which bounds the decay of the Yamabe energy
under a solution to (\ref{flow}).

\begin{lemma} \label{le:yamabedecay} Let $(M^4, g(t))$ be a solution to
(\ref{flow}).  Then
\begin{align*}
\left(\int_M s dV \right)(t) \geq&\ \left( \int_M s dV \right)(0) - C
t^{\frac{1}{2}} \mathcal F(0) \left( \mathcal F(0) - \mathcal F(t) \right).
\end{align*}
\begin{proof} Recall that if $g(t)$ is a one-parameter family of metrics with
$\dt g = h$, then
\begin{align*}
 \dt s =&\ - \gD \tr h + \divg \divg h - \left<h, \Rc \right>.
\end{align*}
Note that $\divg \grad \FF = 0$ as a consequence of
diffeomorphism invariance of $\mathcal \FF$.  Thus for a solution to
(\ref{flow}) we conclude
\begin{align*}
\dt \int_M s =&\ \int_M \left(\left< \grad \FF, \Rc \right> - \frac{1}{2} s \tr
\grad \FF \right) dV
\end{align*}
We directly estimate
\begin{align*}
\brs{\dt \int_M s dV} \leq&\ \brs{ \int_M \left(\left< \grad \FF, \Rc \right> -
\frac{1}{2} s \tr \grad \FF \right) dV}\\
\leq&\ C \brs{\brs{\Rc}}_{L^2} \brs{\brs{\grad \FF}}_{L^2}.
\end{align*}
Thus we may integrate in time to yield
\begin{align*}
\left( \int_M s dV \right)(t) - \left( \int_M s dV \right)(0) \geq&\ - C
\int_0^t \brs{\brs{\Rc}}_{L^2} \brs{\brs{\grad \FF}}_{L^2} dt\\
\geq&\ - C \mathcal F(0) \left( \int_0^t dt \right)^{\frac{1}{2}} \left(
\int_0^t \brs{\brs{\grad \FF}}_{L^2}^2 dt \right)\\
=&\ - C t^{\frac{1}{2}} \mathcal F(0) \left( \mathcal F(0) - \mathcal F(t)
\right).
\end{align*}
\end{proof}
\end{lemma}

So, consider $g = g_i$ some element of the above sequence.  Suppose $T \leq 1$
is
the maximal existence time of the flow $g(t)$.
By the gradient flow property and equation (\ref{eq:energy}) we conclude that
\begin{align*}
 \brs{\brs{z}}_{L^2}^2 (T) \leq \ge_i
\end{align*}
for all $t \in [0, T]$.  Also, applying Lemma \ref{le:yamabedecay} we conclude,
since $T \leq 1$, that
\begin{align*}
\int_M s dV (T) \geq&\ \int_M s dV(0) - C \ge_i\\
\geq&\ Y_{[g(0)]} - C \ge_i.
\end{align*}
Since $Y_{[g(0)]} \geq \sqrt{384 \pi^2 - 96 \ge_i}$ by Proposition \ref{sobprop}, it follows from H\"older's inequality that
\begin{align*}
\left( \sqrt{384 \pi^2} - C \ge_i \right)^2 \leq&\ \left[\left( \int_M s dV \right)(T)
\right]^2\\
\leq&\ \int_M s^2 dV(T).
\end{align*}
It now follows from the Gauss-Bonnet formula that
\begin{align*}
\brs{\brs{W}}_{L^2}^2(T) \leq C \ge_i.
\end{align*}
In particular, we have now shown that there is a universal constant $C$ so that
\begin{align} \label{eq:mainloc10}
\brs{\brs{\ohat{\Rm}}}_{L^2}^2(T) \leq C \ge_i.
\end{align}
Given this, we return to Proposition \ref{sobprop} to conclude that the Sobolev constant is bounded on $[0, T]$.  Suppose
\begin{align*}
\limsup_{t \to T} \brs{\brs{\Rm}}_{\infty} \leq C.
\end{align*}
Since the curvature and Sobolev constants are bounded, it follows from
\cite{Streets} Theorem 6.2 that the flow exists smoothly up to
time $T$, and hence past it, contradicting maximality of $T$.  Therefore we
conclude that
\begin{align*}
\limsup_{t \to T} \brs{\brs{\Rm}}_{\infty} = \infty.
\end{align*}
Let $(x_j, t_j)$ be a sequence of points such that $t_j \to T$ and
\begin{align*}
\limsup_{t \to T} \brs{\brs{\Rm}}_{\infty} = \lim_{j \to \infty}
\brs{\brs{\Rm}}_{\infty} (x_j, t_j) =: \gl_j
\end{align*}
Let
\begin{align*}
g_j(t, x) := \gl_j g\left(t_j + \frac{t}{\gl_j^2}, x \right)
\end{align*}
Consider the sequence of pointed Riemannian manifolds $(M, g_j(t), x_j)$.  They
have uniformly bounded curvatures on the time interval $[- t_j \gl_j^2, 0]$ and
uniformly bounded Sobolev constants, and hence by Theorem 7.1 of \cite{Streets}
we conclude subsequential convergence to a solution $(M_{\infty}, g_{\infty}(t),
x_{\infty})$ of (\ref{flow}) on the time interval $[- \infty, 0]$.  Note that
quadratic curvature functionals are scaling invariant on $M$, so by Fatou's
lemma upper bounds on such integrals pass to the limit $g_{\infty}$.  Moreover,
again using that $\mathcal F(g)$ is scaling invariant, we conclude that
\begin{align*}
\int_{-1}^0 \int_{M_{\infty}} \brs{\grad \FF_{\infty}}^2 dV_{\infty} dt \leq&\
\lim_{j \to \infty} \int_{\frac{-1}{\gl_j^2} + t_j}^{t_j} \int_M \brs{\grad
\FF_j}^2 dV_j dt\\
=&\ \lim_{j \to \infty} \mathcal F\left(g \left(\frac{-1}{\gl_j^2} + t_j\right)
\right) - \mathcal F(g(t_j))\\
=&\ 0.
\end{align*}
Therefore $g_{\infty}(t) = g_{\infty}(0)$ is a critical metric for all $t$ with
$\brs{\brs{\Rm_{\infty}}}_{\infty} = 1$.  Furthermore, we conclude that $\tr
\grad \FF_{\infty} = \gD s = 0$.  Since $\int_M s^2_{\infty} dV_{\infty} \leq
C$, it follows by the maximum principle that $s$ is constant, and this constant
must be zero. Moreover, the limiting manifold is noncompact, and satisfies the
Sobolev inequality
\begin{align*}
\brs{\brs{u}}_{L^4} \leq C \brs{\brs{\N u}}_{L^2}
\end{align*}
where the constant $C$ is bounded uniformly the Sobolev constants of the
metrics $g_i$.  Therefore we may apply Theorem \ref{noncompactgapthm} to
conclude that for $\ge_i$ small enough, $\brs{\brs{\Rm_{\infty}}}_{\infty} = 0$,
a contradiction.  Thus
\begin{align*}
\limsup_{t \to T} \brs{\brs{\Rm}}_{\infty} \leq C
\end{align*}
Thus $T$ is not the maximal existence time, and we have shown that for
sufficiently large $i$ the solution to (\ref{flow}) with initial condition $g_i$
exists at least on $[0, 1]$.
Note that it follows from the above argument that there is a constant $C > 0$ so
that
\begin{align} \label{eq:mainloc20}
 \sup_{M \times \left[\frac{1}{2}, 1 \right]} \brs{\brs{\Rm(g_i)}}_{\infty} \leq
K.
\end{align}
Indeed, if this were not the case, one could choose a sequence $i \to \infty$
and points $(x_i, t_i)$, $t_i \in \left[\frac{1}{2}, 1 \right]$ and repeat the
blowup process.  The resulting blow-up
metric will be critical since $\ge_i \to 0$, and then another application of
Theorem \ref{noncompactgapthm} provides the contradiction.  It is important to
note that we do not have any a priori control over this constant $K$, we merely
know it exists.  Since the curvatures and Sobolev constant are bounded on
$[\frac{1}{2}, 1]$, it follows from \cite{Streets} Theorem 5.4 and the Sobolev
inequality that there exist constants $C_m$ such that
\begin{align} \label{eq:mainloc30}
\sup_{M \times \left[\frac{3}{4}, 1 \right]} \brs{\brs{\N^m \Rm (g_i)}}_{\infty}
\leq C_m C_S K^{m+5}
\end{align}

To finish the first step we show that for $i$ sufficiently large the scalar
curvature of $g_i\left(1 \right)$ is bounded away from zero.  Using the above
estimates, it is clear that if we fix $x \in M$ the sequence of pointed
Riemannian manifolds $\{M, g_i\left(1 \right), x \}$ has a subsequence
which converges, up
to diffeomorphisms, to a new smooth metric $g_{\infty}$.  Since $\ge_i \to 0$,
it follows from the above estimates that $g_{\infty}$ satisfies
$\ohat{\Rm}_{\infty} \equiv 0$, and it then follows from Schur's lemma that
$s_{\infty}$ is constant, and in particular $g_{\infty}$ is isometric to
$g_{S^4}$.  This metric has constant positive
scalar curvature, and since a lower bound on
scalar curvature is diffeomorphism invariant, we conclude that given $\gd > 0$,
for $i$ sufficiently large one has
\begin{align} \label{eq:mainloc40}
 s_{g_i\left(1\right)} \geq s_{g_{S^4}} - \gd.
\end{align}
The second main step is completed in Proposition \ref{analyticconv}, where it is
shown that for $\ge_i$ sufficiently small with respect to $K$, metrics satisfying
(\ref{eq:mainloc10}), (\ref{eq:mainloc20}), (\ref{eq:mainloc30}) and
(\ref{eq:mainloc40}), the solution to (\ref{flow}) exists for all time and
converges exponentially to $g_{S^4}$.  This contradicts the initial hypothesis,
and finishes the proof of the theorem.
\end{proof}

\section{A-priori \texorpdfstring{$L^2$}{L2} growth estimate for
\texorpdfstring{$\grad \FF$}{grad F}}
In this section we give a bound on the growth of $\brs{\brs{\grad \FF}}_{L^2}$
over time intervals of small energy decay.  This is the key input in showing
exponential convergence of long-time solutions of (\ref{flow}) near round
metrics.  To simplify notation we will set
\begin{gather*}
\gf := \grad \FF.
\end{gather*}
The estimate applies in a more general situation which we describe now.  Let
$\ge_0$ be a small constant which will be fixed later and fix a
time interval $[t_0, t_1]$ such that some solution to (\ref{flow}) exists on
$[t_0, t_1]$, has unit volume, and satisfies
\begin{align} \label{eq:gradf5}
\int_{t_0}^{t_1} \int_M \brs{\gf}^2 dV dt \leq \ge \leq 1.
\end{align}
Note that this condition is satisfied for arbitrary time intervals if the
initial condition satisfies
\begin{align*}
\brs{\brs{z}}_{L^2}^2 \leq \ge.
\end{align*}
Furthermore assume that for any $t \in [t_0, t_1]$ one has 
\begin{gather} \label{eq:gradf10}
C_S(g_t) \leq A.
\end{gather}
Without loss of generality we assume $A \geq 1$.  In this setting we derive an
estimate for the $L^2$ norm of $\gf$.  A direct calculation (see \cite{Streets2}
Lemma 13) yields
\begin{gather} \label{L2gradFev}
\begin{split}
\dt \brs{\brs{\gf}}_{L^2}^2 =&\ - \brs{\brs{\gD \gf}}_{L^2}^2 +
\int_M \gf * \N^2 \gf* \Rm\\
&\ + \int_M \gf * \N \gf * \N \Rm + \gf^{*2} * \Rm^{*2} +
\gf^{*2} * \N^2 \Rm.
\end{split}
\end{gather}
Integrating by parts and commuting derivatives yields
\begin{align*}
\brs{\brs{\gD \gf}}_{L^2}^2 =&\ \int_M \N_i \N_i \gf_{j k} \N_l \N_l
\gf_{jk}\\
=&\ \int_M \brs{\N^2 \gf}^2 + \Rm * \N \gf^{*2}\\
\geq&\ \int_M \brs{\N^2 \gf}^2 - C \int_M \brs{\gf} \brs{\N^2 \gf}
\brs{\Rm}\\
&\ - C \int_M \brs{\gf} \brs{\N \gf} \brs{\N \Rm}.
\end{align*}
Combining this with (\ref{L2gradFev}) and integrating over the time interval
$[t_0, t_1]$ yields
\begin{gather} \label{semainbound}
\begin{split}
\brs{\brs{\gf}}^2_{L^2(g_{t_1})} & + \int_{t_0}^{t_1} \int_M \brs{\N^2
\gf}^2\\
\leq&\ \brs{\brs{\gf}}^2_{L^2(g_{t_0})} + C \int_{t_0}^{t_1} \int_M \left[
\brs{\gf} \brs{\N^2 \gf} \brs{\Rm} \right.\\
&\ \left. + \brs{\gf} \brs{\N \gf} \brs{\N \Rm} + \brs{\gf}^2
\brs{\Rm}^2 + \brs{\gf}^2 \brs{\N^2 \Rm} \right].
\end{split}
\end{gather}
\noindent We now proceed to bound the terms on the right hand side of the above
inequality in a series of lemmas.
\begin{lemma} \label{lemmabound1} Given $(M^4, g(t))$ a solution to (\ref{flow})
satisfying (\ref{eq:gradf5}) and (\ref{eq:gradf10}), there is a constant $C$
depending on $\mathcal F(g(t_1))$ such that 
\begin{align*}
\int_{t_0}^{t_1} \int_M \brs{\gf}^2 \brs{\Rm}^2 \leq C A^2 \ge^{\frac{1}{2}}
\left[ 1 + \sup_{t_0 \leq
t \leq t_1} \int_M \brs{\N^2 \Rm}^2 +
 \int_{t_0}^{t_1} \int_M \brs{\N^2 \gf}^2 \right].
\end{align*}
\begin{proof} First we apply H\"older's inequality and the Sobolev inequality to
yield
\begin{gather} \label{sebound1}
\begin{split}
\int_M \brs{\gf}^{2} \brs{\Rm}^2 \leq&\ \left( \int_M \brs{\gf}^4
\right)^{\frac{1}{2}} \left( \int_M \brs{\Rm}^4 \right)^{\frac{1}{2}}\\
\leq&\ C A^2 \left( \int_M \brs{\N \gf}^2 + \int_M \brs{\gf}^2 \right)
\left( \int_M \brs{\N \Rm}^2 + \int_M \brs{\Rm}^2 \right)\\
=&\ I + II + III + IV, 
\end{split}
\end{gather}
where the Roman numerals refer to the four different terms in the expanded
product.  First we bound the main term, integrating by parts
\begin{align*}
I =&\ C A^2 \int_M \brs{\N \gf}^2 \int_M \brs{\N \Rm}^2\\
=&\ C A^2 \left( \int_M \left< \gf, \gD \gf \right> \right) \left(
\int_M \left< \Rm, \gD \Rm \right> \right)\\
\leq&\ C A^2 \left( \int_M \brs{\gf}^2 \right)^{\frac{1}{2}} \left( \int_M
\brs{\gD \gf}^2 \right)^{\frac{1}{2}} \left( \int_M \brs{\Rm}^2
\right)^{\frac{1}{2}} \left( \int_M \brs{\gD \Rm}^2 \right)^{\frac{1}{2}}\\
\leq&\ C A^2 \left( \int_M \brs{\gf}^2 \right)^{\frac{1}{2}} \left( \int_M
\brs{\N^2 \gf}^2 \right)^{\frac{1}{2}} \left( \int_M \brs{\N^2 \Rm}^2
\right)^{\frac{1}{2}}.
\end{align*}
Integrating this bound in time and applying H\"olders inequality to the time
integral yields
\begin{gather} \label{sebound5}
\begin{split}
\int_{t_0}^{t_1} I \leq&\ C A^2 \sup_{t_0 \leq
t \leq t_1} \left( \int_M \brs{\N^2 \Rm}^2 \right)^{\frac{1}{2}}
\left(\int_{t_0}^{t_1} \int_M \brs{\gf}^2 \right)^{\frac{1}{2}} \left(
\int_{t_0}^{t_1} \int_M \brs{\N^2 \gf}^2 \right)^{\frac{1}{2}}\\
\leq&\ C A^2 \ge^{\frac{1}{2}} \left[ \sup_{t_0 \leq
t \leq t_1} \int_M \brs{\N^2 \Rm}^2  + 
 \int_{t_0}^{t_1} \int_M \brs{\N^2 \gf}^2 \right].
\end{split}
\end{gather}
The lower order terms are easier to bound.  We bound by interpolation
\begin{align*}
\int_{t_0}^{t_1} II =&\ \int_{t_0}^{t_1} \int_M \brs{\N \gf}^2 \int_M
\brs{\Rm}^2\\
\leq&\ C \int_{t_0}^{t_1} \left( \int_M \brs{\gf}^2 \right)^{\frac{1}{2}}
\left( \int_M \brs{\N^2 \gf}^2 \right)^{\frac{1}{2}}\\
\leq&\ C \left(\int_{t_0}^{t_1} \int_M \brs{\gf}^2 \right)^{\frac{1}{2}}
\left( \int_{t_0}^{t_1} \int_M \brs{\N^2 \gf}^2 \right)^{\frac{1}{2}}\\
\leq&\ C \ge^{\frac{1}{2}} \left[ 1 + \int_{t_0}^{t_1} \int_M \brs{\N^2 E}^2
\right].
\end{align*}
For the third term we again interpolate
\begin{align*}
\int_{t_0}^{t_1} III =&\ C \int_{t_0}^{t_1} \int_M \brs{\gf}^2 \int_M \brs{\N
\Rm}^2\\
\leq&\ C \int_{t_0}^{t_1} \int_M \brs{\gf}^2 \left( \int_M \brs{\Rm}^2
\right)^{\frac{1}{2}} \left(\int_M \brs{\N^2 \Rm}^2 \right)^{\frac{1}{2}}\\
\leq&\ C \sup_{t_0 \leq t \leq t_1} \left(\int_M \brs{\N^2
\Rm}^2 \right)^{\frac{1}{2}} \int_{t_0}^{t_1} \int_M \brs{\gf}^2\\
\leq&\ C \ge^{\frac{1}{2}} \left[ 1 + \sup_{t_0 \leq t \leq t_1} \int_M
\brs{\N^2 \Rm}^2 \right].
\end{align*}
Finally we make the bound
\begin{align*}
\int_{t_0}^{t_1} IV =&\ C \int_{t_0}^{t_1} \int_M \brs{\gf}^2 \int_M
\brs{\Rm}^2\\
\leq&\ C \int_{t_0}^{t_1} \int_M \brs{\gf}^2\\
\leq&\ C \ge.
\end{align*}
Combining these bounds gives the result.
\end{proof}
\end{lemma}

\begin{lemma} \label{lemmabound2} Given $(M^4, g(t))$ a solution to (\ref{flow})
satisfying (\ref{eq:gradf5}) and (\ref{eq:gradf10}), there is a constant $C$
depending on $\mathcal F(g(t_1))$ such that 
\begin{gather} \label{sebound10}
\begin{split}
\int_{t_0}^{t_1} \int_M \brs{\gf}^2 \brs{\N^2 \Rm} \leq C A^2
\ge^{\frac{1}{2}} \left[1 +  \sup_{t_0 \leq
t \leq t_1} \int_M \brs{\N^2 \Rm}^2 + \int_{t_0}^{t_1} \int_M \brs{\N^2 \gf}^2
\right].
\end{split}
\end{gather}
\begin{proof} First we apply H\"older's inequality and the Sobolev inequality to
bound
\begin{align*}
\int_M \brs{\gf}^{2} \brs{\N^2 \Rm} \leq&\ \int_M \brs{\gf}^2
\brs{\N^2
\Rm}\\
\leq&\ \left( \int_M \brs{\gf}^4 \right)^{\frac{1}{2}} \left( \int_M
\brs{\N^2 \Rm}^2 \right)^{\frac{1}{2}}\\
\leq&\ A \left( \int_M \brs{\N \gf}^2 + \int_M \brs{E}^2 \right) \left( \int_M
\brs{\N^2
\Rm}^2 \right)^{\frac{1}{2}}\\
\leq&\ C A \left[ \left( \int_M \brs{\gf}^2 \right)^{\frac{1}{2}} \left( \int_M
\brs{\N^2 \gf}^2 \right)^{\frac{1}{2}} + \int_M \brs{\gf}^2 \right] \left(
\int_M \brs{\N^2 \Rm}^2
\right)^{\frac{1}{2}}.
\end{align*}
In the last line we applied interpolation to the integral $\int_M \brs{\N
\gf}^2$.  The second term above may be integrated in time to yield
\begin{align*}
C A \int_{t_0}^{t_1} \int_M \brs{\gf}^2 \left(\int_M \brs{\N^2 \Rm}^2
\right)^{\frac{1}{2}} \leq&\ C A \ge \left(\sup_{t_0 \leq t \leq t_1} \int_M
\brs{\N^2 \Rm}^2 \right)^{\frac{1}{2}}\\
\leq&\ C A \ge \left[1 + \sup_{t_0 \leq t \leq t_1} \int_M \brs{\N^2 \Rm}^2
\right].
\end{align*}
The first term above is integrated in time and bounded as in line
(\ref{sebound5}), yielding the result.
\end{proof}
\end{lemma}

\begin{lemma} \label{lemmabound3} Given $(M^4, g(t))$ a solution to (\ref{flow})
satisfying (\ref{eq:gradf5}) and (\ref{eq:gradf10}), there is a constant $C$
depending on $\mathcal F(g(t_1))$ such that 
\begin{align*}
\int_{t_0}^{t_1} \int_M \brs{E} \brs{\N E} \brs{\N \Rm} \leq C A^2
\ge^{\frac{1}{2}} \left[1 +  \sup_{t_0 \leq
t \leq t_1} \int_M \brs{\N^2 \Rm}^2 + \int_{t_0}^{t_1} \int_M \brs{\N^2 \gf}^2
\right].
\end{align*}
\begin{proof}  We apply H\"older's inequality, the Sobolev inequality and
interpolation to bound
\begin{align*}
\int_M \brs{\gf} \brs{\N \gf} \brs{\N \Rm} \leq&\ \left(\int_M
\brs{\gf}^2
\right)^{\frac{1}{2}} \left( \int_M \brs{\N \gf}^4 \right)^{\frac{1}{4}}
\left( \int_M \brs{\N \Rm}^4 \right)^{\frac{1}{4}}\\
\leq&\ C A^2 \left( \int_M \brs{\gf}^2 \right)^{\frac{1}{2}} \left( \int_M
\brs{\N^2 \gf}^2 + \int_M \brs{\N E}^2 \right)^{\frac{1}{2}}\cdot\\
&\  \left( \int_M \brs{\N^2 \Rm}^2
 + \int_M \brs{\N \Rm}^2 \right)^{\frac{1}{2}}\\
\leq&\ C A^2 \left( \int_M \brs{\gf}^2 \right)^{\frac{1}{2}} \left( \int_M
\brs{\N^2 \gf}^2 + \int_M \brs{E}^2 \right)^{\frac{1}{2}}\cdot\\
&\  \left( \int_M \brs{\N^2 \Rm}^2
 + \int_M \brs{\Rm}^2 \right)^{\frac{1}{2}}.
\end{align*}
The time integral of each of the terms above has been bounded in the previous
two lemmas, and so the result follows.
\end{proof}
\end{lemma}

\begin{lemma} \label{lemmabound4} Given $(M^4, g(t))$ a solution to (\ref{flow})
satisfying (\ref{eq:gradf5}) and (\ref{eq:gradf10}), there is a constant $C$
depending on $\mathcal F(g(t_1))$ such that 
\begin{align*}
\int_{t_0}^{t_1} \int_M \brs{E} \brs{\N^2 E} \brs{\Rm} \leq&\ C A^2
\ge^{\frac{1}{4}} \left[ 1 + \sup_{t_0 \leq t \leq t_1} \int_M \brs{\N^2
\Rm}^2 + \int_{t_0}^{t_1} \int_M \brs{\N^2 E}^2 \right].
\end{align*}
\begin{proof} We start by applying H\"older's inequality and the Sobolev
inequality to bound
\begin{align*}
\int_M \brs{E} \brs{\N^2 E} \brs{\Rm} \leq&\ \left(\int_M \brs{E}^4
\right)^{\frac{1}{4}} \left( \int_M \brs{\N^2 E}^2 \right)^{\frac{1}{2}} \left(
\int_M \brs{\Rm}^4 \right)^{\frac{1}{4}}\\
\leq&\ A^2 \left( \int_M \brs{\N E}^2 + \int_M \brs{E}^2 \right)^{\frac{1}{2}}
\left( \int_M \brs{\N^2 E}^2 \right)^{\frac{1}{2}} \cdot\\
&\ \left( \int_M \brs{\N \Rm}^2 + \int_M \brs{\Rm}^2 \right)^{\frac{1}{2}}\\
=&\ A^2 \left(I + II + III + IV \right).
\end{align*}
where the Roman numerals denote the four terms in the expanded product above
after applying the inequality $\sqrt{a + b} \leq \sqrt{a} + \sqrt{b}$.  First we
bound the highest order term
\begin{align*}
\int_{t_0}^{t_1} I =&\ \int_{t_0}^{t_1} \left(\int_M \brs{\N E}^2
\right)^{\frac{1}{2}} \left( \int_M \brs{\N^2 E}^2 \right)^{\frac{1}{2}} \left(
\int_M \brs{\N \Rm}^2 \right)^{\frac{1}{2}}\\
\leq&\ \int_{t_0}^{t_1} \left( \int_M \brs{E}^2 \right)^{\frac{1}{4}} \left(
\int_M \brs{\N^2 E}^2 \right)^{\frac{3}{4}} \left( \int_M \brs{\Rm}^2
\right)^{\frac{1}{4}} \left( \int_M \brs{\N^2 \Rm}^2 \right)^{\frac{1}{4}}\\
\leq&\ C \sup_{t_0 \leq t \leq t_1} \left( \int_M \brs{\N^2 \Rm}^2
\right)^{\frac{1}{4}} \left(\int_{t_0}^{t_1} \int_M \brs{E}^2
\right)^{\frac{1}{4}} \left( \int_{t_0}^{t_1} \int_M \brs{\N^2 E}^2
\right)^{\frac{3}{4}}\\
\leq&\ C \ge^{\frac{1}{4}} \left[ \sup_{t_0 \leq t \leq t_1} \int_M \brs{\N^2
\Rm}^2 + \int_{t_0}^{t_1} \int_M \brs{\N^2 E}^2 \right].
\end{align*}
Next we bound
\begin{align*}
\int_{t_0}^{t_1} II =&\ \int_{t_0}^{t_1} \left(\int_M \brs{\N
E}^2\right)^{\frac{1}{2}} \left(\int_M \brs{\Rm}^2 \right)^{\frac{1}{2}} \left(
\int_M \brs{\N^2 E}^2 \right)^{\frac{1}{2}}\\
\leq&\ C \int_{t_0}^{t_1} \left( \int_M \brs{E}^2 \right)^{\frac{1}{4}} \left(
\int_M \brs{\N^2 E}^2 \right)^{\frac{3}{4}}\\
\leq&\ C \left( \int_{t_0}^{t_1} \int_M \brs{E}^2 \right)^{\frac{1}{4}} \left(
\int_{t_0}^{t_1} \int_M \brs{\N^2 E}^2 \right)^{\frac{3}{4}}\\
\leq&\ C \ge^{\frac{1}{4}} \left[1 + \int_M \brs{\N^2 E}^2 \right].
\end{align*}
For the third term we bound
\begin{align*}
\int_{t_0}^{t_1} III =&\ \int_{t_0}^{t_1} \left(\int_M \brs{E}^2
\right)^{\frac{1}{2}} \left( \int_M \brs{\N^2 E}^2 \right)^{\frac{1}{2}} \left(
\int_M \brs{\N \Rm}^2 \right)^{\frac{1}{2}}\\
\leq&\ \int_{t_0}^{t_1} \left(\int_M \brs{E}^2 \right)^{\frac{1}{2}} \left(
\int_M \brs{\N^2 E}^2 \right)^{\frac{1}{2}} \left( \int_M \brs{\N^2 \Rm}^2
\right)^{\frac{1}{4}} \left( \int_M \brs{\Rm}^2 \right)^{\frac{1}{4}}\\
\leq&\ C \left(\sup_{t_0 \leq t \leq t_1} \int_M \brs{\N^2 \Rm}^2
\right)^{\frac{1}{4}} \left(\int_{t_0}^{t_1} \int_M \brs{E}^2
\right)^{\frac{1}{2}} \left( \int_{t_0}^{t_1} \int_M \brs{\N^2 E}^2
\right)^{\frac{1}{2}}\\
\leq&\ C \ge^{\frac{1}{2}} \left[ 1 + \sup_{t_0 \leq t \leq t_1} \int_M
\brs{\N^2 \Rm}^2 + \int_{t_0}^{t_1} \int_M \brs{\N^2 E}^2 \right]. 
\end{align*}
Finally we estimate
\begin{align*}
\int_{t_0}^{t_1} IV =&\ \int_{t_0}^{t_1} \left(\int_M \brs{E}^2
\right)^{\frac{1}{2}} \left( \int_M \brs{\N^2 E}^2 \right)^{\frac{1}{2}} \left(
\int_M \brs{\Rm}^2 \right)^{\frac{1}{2}}\\
\leq&\ C \left( \int_{t_0}^{t_1} \int_M \brs{E}^2 \right)^{\frac{1}{2}} \left(
\int_{t_0}^{t_1} \int_M \brs{\N^2 E}^2 \right)^{\frac{1}{2}}\\
\leq&\ C \ge^{\frac{1}{2}} \left[1 + \int_{t_0}^{t_1} \int_M \brs{\N^2 E}^2
\right].
\end{align*}
Combining these four estimates and using that $\ge \leq 1$ gives the result.
\end{proof}
\end{lemma}

\begin{prop} \label{gradFbound}  Given $(M^4, g(t))$ a solution to (\ref{flow})
satisfying (\ref{eq:gradf5}) and (\ref{eq:gradf10}), there is a constant $C > 0$
depending on $\mathcal F(g(t_1))$ so that if $\ge$ is chosen small with
respect to $A$ and $\mathcal F(g(t_1))$ one has
\begin{align*}
\sup_{t_0 \leq t \leq t_1} \brs{\brs{\gf}}^2_{L^2} + \int_{t_0}^{t_1}
\int_M \brs{\N^2 \gf}^2 \leq&\ 2 \brs{\brs{\gf}}_{L^2(g_{t_0})}^2 + C A^2
\ge^{\frac{1}{4}} \left[1 + \sup_{t_0 \leq t \leq t_1} \int_M \brs{\N^2 \Rm}^2
\right].
\end{align*}
\begin{proof}
Combining Lemmas \ref{lemmabound1} - \ref{lemmabound4} and plugging into
(\ref{L2gradFev}) yields
\begin{align*}
\sup_{t_0 \leq t \leq t_1} \brs{\brs{\gf}}^2_{L^2}& + \int_{t_0}^{t_1}
\int_M \brs{\N^2 \gf}^2\\
\leq&\ \brs{\brs{\gf}}_{L^2(g_{t_0})}^2 + C A^2 \ge^{\frac{1}{4}} \left[ 1 +
\sup_{t_0
\leq t \leq
t_1} \int_M \brs{\N^2 \Rm}^2 + 
\int_{t_0}^{t_1} \int_M \brs{\N^2 \gf}^2 \right]
\end{align*}
Therefore for $\ge$ chosen small enough with respect to $A$ and the constants of
the lemmas, which depend on $\mathcal F(g(t_1))$, we conclude the result.
\end{proof}
\end{prop}

\section{Exponential Convergence}

\begin{prop} \label{analyticconv} Given $K > 0$, $0 < \gd << 1$, there exists
$\ge > 0$
so that if $(M^4, g(t))$ is a solution to (\ref{flow}) which exists on $[0, 1]$
and satisfies $Y_{[g(0)]} > 0$,
\begin{gather} \label{eq:anal10}
 \begin{split}
\sup_{t \in [0,1]} \brs{\brs{\ohat{\Rm}}}_{L^2}^2(g(t)) \leq&\ \ge,\\
 \end{split}
\end{gather}
(\ref{eq:mainloc20}), (\ref{eq:mainloc30}),
and (\ref{eq:mainloc40}), then the solution exists for all time and converges to
either $g_{S^4}$ or $g_{\mathbb R \mathbb P^4}$.
\begin{proof}  The strategy is to use the key coercivity estimate of
Proposition \ref{maincoerciveestimate} to show exponential decay of
$\grad \FF$. 
With this decay in hand, an argument exploiting a multiplicative Sobolev
inequality and Moser iteration
can be applied to conclude exponential convergence of the flow.

Let $(M^4, g(t))$ be a solution to (\ref{flow}) satisfying the hypotheses of the
proposition.  As in section 4, assume by passing to the double cover that $M$ is oriented.  Observe that $g(1)$ trivially satisfies by hypothesis
\begin{gather} \label{eq:main10}
\begin{split}
s >&\ s_{g_{S^4}} - 2 \gd\\
\brs{\brs{\Rm}}_{\infty} <&\ 2 K\\
\brs{\brs{\ohat{\Rm}}}_{L^2}^2 <&\ 2 \ge
\end{split}
\end{gather}
Let $\Omega = \{ t \in [1, \infty) | (\ref{eq:main10}) \mbox{ is satisfied} \}$.
 $\Omega$ is certainly open, and we aim to show that $\Omega$ is closed.  Let $T
\in \Omega$.  If $\ge$ is small enough, we may apply Proposition \ref{sobprop}
to conclude that there is a uniform constant $A$ such that
\begin{align*}
\sup_{t \in [0, T]} C_S(g(t)) \leq A.
\end{align*}
Likewise, using (\ref{eq:mainloc20}), we have that $\sup_{M \times [T -
\frac{1}{2}, T]} \brs{\brs{\Rm}}_{\infty} \leq 2K$.  Using this and the
curvature bound of (\ref{eq:main10}), we may argue as in section 4 using the
derivative estimates for solutions to (\ref{flow}) to conclude that there are
constants $C_m$ such that
\begin{align} \label{eq:main12}
 \brs{\brs{\N^m \Rm(g(T))}}_{\infty} \leq C_m C_S K^{m+5}.
\end{align}
Thus if condition (\ref{eq:main10}) holds on $[0, T)$, the solution to
(\ref{flow}) exists smoothly up to and past time $T$.

We now derive exponential decay of $\brs{\brs{\grad \FF}}_{L^2}$.  First note
that, using (\ref{eq:energy}), we have that
\begin{align*}
 \dt \brs{\brs{z}}_{L^2}^2 =&\ \dt \left( 8 \pi^2 \chi(M) +
\brs{\brs{z}}_{L^2}^2 \right) = \dt \mathcal F = - \brs{\brs{\grad
\FF}}_{L^2}^2.
\end{align*}
Provided say $\gd < \frac{1}{100}$, by Proposition \ref{maincoerciveestimate} we
conclude
that there is a constant
$\eta > 0$ so that if $\ge$ is chosen small with respect to $A$ and $K$, then
given $t \geq 1 \in \Omega$, we have
\begin{align*}
 \brs{\brs{\grad \FF}}_{L^2}^2 \geq \eta \brs{\brs{z}}_{L^2}^2.
\end{align*}
Combining this with the line above we conclude that for $t \geq 1,$
\begin{align}\label{eq:main45}
\brs{\brs{z}}_{L^2}^2(t) \leq \ge e^{-\eta t}.
\end{align}
Given exponential decay of the energy, it is natural to expect exponential decay
of its time derivative.  We claim that there exists a constant $P = P(A, K)$
such that for $t \geq 1$,
\begin{align} \label{eq:main50}
\brs{\brs{\grad \FF}}_{L^2}^2(t) <&\ P \ge^{\frac{1}{4}} e^{-\frac{\eta}{4} t}.
\end{align}
We first need to show this estimate on the time interval $[1, \frac{5}{4}]$. 
Note that
\begin{align*}
\int_{\frac{3}{4}}^1 \brs{\brs{\grad \FF}}_{L^2}^2 = \mathcal F(\frac{3}{4}) -
\mathcal F(1) \leq \ge.
\end{align*}
Thus there exists $s, \frac{3}{4} \leq s \leq 1$ such that $\brs{\brs{\grad
\FF}}_{L^2}^2(s) \leq \ge$.  Using Proposition \ref{gradFbound} and
(\ref{eq:mainloc30}) we conclude that if $\ge$ is chosen small enough we have
\begin{align*}
\sup_{t \in [1, \frac{5}{4}]} \brs{\brs{\grad \FF}}_{L^2}^2 \leq C
\ge^{\frac{1}{4}}
\end{align*}
which proves (\ref{eq:main50}) on $[1, \frac{5}{4}]$ for $P = C e^{\frac{5
\eta}{15}}$.  Next we show (\ref{eq:main50}) for arbitrary times $t \geq
\frac{5}{4}$.
Observe for any $1 \leq t_1 \leq t_2$ the estimate
\begin{gather} \label{eq:main55}
\begin{split}
\int_{t_1}^{t_2} \brs{\brs{\grad \FF}}_{L^2}^2 =&\ \mathcal F(t_1) - \mathcal
F(t_2)\\
=&\ \left( 8 \pi^2 \chi(M) + \brs{\brs{z}}_{L^2}^2(t_1) \right) - \left( 8 \pi^2
\chi(M) + \brs{\brs{z}}_{L^2}^2(t_2) \right)\\
\leq&\ \brs{\brs{z}}_{L^2}^2(t_1)\\
\leq&\ \ge e^{-\eta t_1}.
\end{split}
\end{gather}
Now fix some $t \geq \frac{5}{4}$.  Applying (\ref{eq:main55}) for $t_1 = t -
\frac{1}{4}, t_2 = t$
we conclude that there exists $s \in \left[t - \frac{1}{4}, t \right]$ such that
$\brs{\brs{\grad
\FF}}_{L^2}^2(s) \leq \ge e^{- \eta\left(t - \frac{1}{4} \right)}$.  Next we
apply Proposition
\ref{gradFbound} with $t_0 = s, t_1 = t$ and apply (\ref{eq:main12}), to
conclude
\begin{align*}
\brs{\brs{\grad \FF}}_{L^2}^2(t) \leq&\ 2 \brs{\brs{\grad \FF}}_{L^2}^2(s) +
C(K) A^2 \ge^{\frac{1}{4}} e^{- \frac{\eta\left(t-\frac{1}{4} \right)}{4}}\\
\leq&\ C(A, K) \ge^{\frac{1}{4}} \left(e^{-\eta\left(t-\frac{1}{4} \right)} +
e^{-\frac{\eta\left(t-\frac{1}{4} \right)}{4}}
\right)\\
\leq&\ C(A, K)\ge^{\frac{1}{4}} e^{-\frac{\eta}{4} t}.
\end{align*}
This finishes the proof of (\ref{eq:main50}).  We now use this estimate to show
that $[0, \infty) \subset \Omega$.  Apply Theorem \ref{multsob} with $p = 8, m =
2$ and $\alpha =
\frac{4}{5}$ and use (\ref{eq:main12}) to estimate
\begin{gather} \label{eq:main35}
\begin{split}
\int_{1}^T \brs{\brs{\grad \FF}}_{\infty} \leq&\ C A \int_{1}^T
\brs{\brs{\grad \FF}}_{L^2}^{\frac{1}{5}} \left( \brs{\brs{\N \grad \FF}}_{L^8}
+ \brs{\brs{\grad \FF}}_{L^8} \right)^{\frac{4}{5}}\\
\leq&\ C(A, K) \ge^{\frac{1}{40}} \int_1^T e^{-\frac{\eta}{40} t}\\
\leq&\ C(A, K) \ge^{\frac{1}{40}}.
\end{split}
\end{gather}
Likewise another application of Theorem \ref{multsob} yields
\begin{align*}
\int_{1}^T \brs{\brs{\N^2 \grad \FF}}_{\infty} \leq C A \int_{1}^T
\brs{\brs{\N^2 \grad \FF}}_{L^2}^{\frac{1}{5}} \left( \brs{\brs{\N^3 \grad
\FF}}_{L^8} + \brs{\brs{\grad \FF}}_{L^8} \right)^{\frac{4}{5}}
\end{align*}
Integrating by parts and applying H\"older's inequality and (\ref{eq:main12}) we
conclude that for $t \geq 1$
\begin{align*}
 \brs{\brs{\N^2 \grad \FF}}_{L^2} \leq&\ \brs{\brs{\grad
\FF}}_{L^2}^{\frac{1}{2}} \brs{\brs{\N^4 \grad \FF}}_{L^2}^{\frac{1}{2}} \leq
C(K)
\brs{\brs{\grad \FF}}_{L^2}^{\frac{1}{2}}.
\end{align*}
Thus we conclude
\begin{gather} \label{eq:main40}
\begin{split} 
\int_{1}^T \brs{\brs{\N^2 \grad \FF}}_{\infty} \leq&\ C(A, K) \int_{1}^T
\brs{\brs{\grad \FF}}_{L^2}^{\frac{1}{10}}\\
\leq&\ C(A, K) \ge^{\frac{1}{80}} \left( \int_{1}^T e^{-\frac{\eta}{80} t}
\right)\\
\leq&\ C(A, K) \ge^{\frac{1}{80}}
\end{split}
\end{gather}
Using these two estimates we can finish the proof.  Recall the evolution
equation computed above,
\begin{align*}
 \dt s =&\ - \gD^2 s - \left<r, \grad \FF \right>.
\end{align*}
Therefore, for times $t \in \Omega$, we conclude using
(\ref{eq:main35}) and (\ref{eq:main40}), for any $x \in M$,
\begin{align*}
s(x,t) - s(x, \tau) \geq&\ - \int_{\tau}^t \brs{\brs{\N^2 \grad \FF}}_{\infty} +
B \brs{\brs{\grad \FF}}_{\infty}\\
\geq&\ - C(A, K) \ge^{\frac{1}{80}}.
\end{align*}
It follows that if $\ge$ is chosen initially small enough, then we may conclude
$s > s_{g_{S^4}} - 2 \gd$ for all times $t \leq T$.  A completely analogous
argument shows that
\begin{align*}
\brs{\brs{\Rm}}_{\infty}(T) \leq \brs{\brs{\Rm}}_{\infty}(1) + C(A, K)
\ge^{\frac{1}{80}}.
\end{align*}
Thus again for $\ge$ chosen small with respect to $A$ and $K$ we conclude
\begin{align*}
\brs{\brs{\Rm}}_{\infty}(T) < 2 K 
\end{align*}
The final bound of (\ref{eq:main10}) follows in an analogous fashion.  Since $T$
was arbitrary, we conclude $[0, \infty) \subset \Omega$.  The
estimates we have shown already imply uniform $C^k$ convergence $g(t) \to
g_{\infty}$ for any $k$.  The decay estimate (\ref{eq:main50}) and the bound $s
> s_{g_{S^4}} - 2 \gd$ together imply that $g_{\infty}$ is a critical metric
with small energy and positive Yamabe constant, which is isometric to $(S^4,
g_{S^4})$ by Theorem \ref{gapthm}.  The proposition follows.
\end{proof}
\end{prop}

\section{Related Questions}

It is tempting to ask what the optimal value of $\ge$ is in the statement of the
three main theorems.  At least for Theorem \ref{mainthm}, it seems natural,
given the main theorem of \cite{CGY}, that $16 \pi^2$ is the optimal value. 
However, this is not completely clear, since solutions to (\ref{flow}) do not
necessarily preserve upper bounds on the Weyl tensor.  Indeed, it was exactly
this problem which forced us to use Lemma \ref{le:yamabedecay} to ensure that
the $L^2$ norm of the Weyl curvature was staying small for a fixed time.

However, if instead of (\ref{flow}), one considered the \emph{Bach flow}, i.e.
the negative gradient flow of the squared $L^2$ norm of the Weyl curvature, then
the hypothesis $\brs{\brs{W}}_{L^2}^2 < 16 \pi^2 \chi(M)$ becomes quite
natural. It is furthermore natural to conjecture in this setting that solutions
to the Bach flow with initial condition satisfying this hypothesis exist for all
time and converge to round metrics.  Many of the techniques used here can likely
be adapted to this setting, but new challenges will certainly arise.  Indeed, to
even define the Bach flow requires adding a certain conformal term to the flow
to overcome the nonparabolicity of the Bach flow which arises due to the
conformal invariance of the Bach tensor.  The existence of this flow with small
energy remains an interesting open question.

\section{Appendix: Sobolev Inequalities}

In this appendix we record a multiplicative Sobolev inequality for Riemannian
manifolds.  The proof is as adaptation of techniques used in \cite{Lady}.

\begin{thm} \label{multsob} Let $(M^4, g)$ be a Riemannian manifold of unit
volume.  For $u \in
C_0^1(M)$, $4 < p \leq \infty, 0 \leq m \leq \infty$ we have
\begin{gather} \label{sob4}
\nm{u}{\infty} \leq C_S \cdot C(n,m,p) \nm{u}{m}^{1 - \ga} \left( \nm{\N u}{p} +
\nm{u}{p} \right)^\ga
\end{gather}
where
$0 < \ga \leq 1$ satisfies $\frac{1}{\ga} = \left( \frac{1}{4} - \frac{1}{p}
\right)m + 1$
\begin{proof} Let $A$ denote the Sobolev constant of $(M, g)$.  Fix $p > 4$, and
rescale $u$ such that
\begin{align*}
A \left( \brs{\brs{\N u}}_{L^p} + \brs{\brs{u}}_{L^p} \right) = 1.
\end{align*}
Let $q = \frac{2p}{p-2}$ and note that for any $w \geq 0,$
\begin{align*}
\brs{\brs{u^{1 + w}}}_{L^4} \leq&\ A \left( \brs{\brs{\N (u^{1+w})}}_{L^2} +
\brs{\brs{u^{1+w}}}_{L^2} \right)\\
\leq&\ A (1 + w) \brs{\brs{u^w}}_{L^q} \left(\brs{\brs{\N u}}_{L^p} +
\brs{\brs{u}}_{L^p} \right)\\
\leq&\ (1+w) \brs{\brs{u^w}}_{L^q}.
\end{align*}
Let $j = \frac{4}{q} \in (2, 4]$.  Then we can rewrite the above estimate as
\begin{align*}
\brs{\brs{u}}_{j(1+w)q} \leq&\ (1 + w)^{\frac{1}{1+w}} \brs{\brs{u}}_{w
q}^{\frac{w}{w+1}}
\end{align*}
We want to apply this estimate inductively.  To that end let $w_0 =
\frac{m}{q}$, $w_{i+1} = j(1 + w_i)$, $\gd_i = \frac{w_i}{w_i + 1}$, $C_i = (1 +
w_i)^{\frac{1}{1 + w_i}}$.  Using this notation the above estimate reads
\begin{align*}
\brs{\brs{u}}_{w_{i+1} q} \leq C_i \brs{\brs{u}}_{w_i q}^{\gd_i}.
\end{align*}
Applying this estimate inductively yields
\begin{align*}
\brs{\brs{u}}_{w_i q} \leq \left( \prod_{l = 0}^{i-1} C_l^{\gd_{l+1} \dots
\gd_{i-1}} \right) \brs{\brs{u}}_m^{\gd_0 \dots \gd_{i-1}}.
\end{align*}
Now observe the formula
\begin{align*}
1 + w_i = j^i w_0 + \sum_{l = 0}^i j^l
\end{align*}
This implies that there exists a constant $C$ depending on $m$ and $p$ such that
\begin{align*}
\frac{1}{C} j^i \leq 1 + w_i \leq C j^i
\end{align*}
Since each $\gd_i \leq 1$ this implies the estimate
\begin{align*}
\log \prod_{l = 0}^{i} C_l^{\gd_{l+1} \dots \gd_i} \leq&\ \sum_{l = 0}^i
\frac{1}{1 + w_l} \log (1 + w_l)\\
\leq&\ \sum_{l = 0}^{\infty} C j^{-l} \left( l \log j \right)\\
\leq&\ C.
\end{align*}
Furthermore we compute
\begin{align*}
\prod_{l = 0}^{\infty} \gd_l =&\ \lim_{i \to \infty} j^i \frac{w_0}{1 + w_i}\\
=&\ \frac{w_0}{w_0 + \frac{j}{j-1}}\\
=&\ 1 - \alpha.
\end{align*}
\end{proof}
\end{thm}

\bibliographystyle{hamsplain}

\end{document}